\documentclass[%
 aip,
 jcp,
 amsmath,amssymb,
 reprint,%
 longbibliography
]{revtex4-1}

\usepackage{graphicx}
\usepackage{dcolumn}
\usepackage{bm}
\usepackage[thinc]{esdiff}
\usepackage{subcaption}
\usepackage{amsmath,amsfonts}
\usepackage{physics}
\usepackage{verbatim}
\usepackage{hyperref}
\hypersetup{colorlinks=true,linkcolor=red,citecolor=blue}
\usepackage{xspace}
\usepackage[dvipsnames]{xcolor}
\usepackage[normalem]{ulem}
\usepackage{multirow}
\usepackage{algorithm,algpseudocode}

\newcommand{\nel}{N_{\mathrm{el}}}
\newcommand{\npw}{N_{\mathrm{pw}}}

\newcommand{\editor}[2]{%
  \expandafter\newcommand\csname #1note\endcsname[1]{%
    \textcolor{#2}{(\textbf{#1:} ##1)}}%
  \expandafter\newcommand\csname #1\endcsname[1]{%
    \textcolor{#2}{##1}}%
  \expandafter\newcommand\csname #1cancel\endcsname[1]{%
    \textcolor{#2}{\sout{##1}}}%
  \expandafter\newcommand\csname #1change\endcsname[2]{%
    \textcolor{#2}{\sout{##1} ##2}}%
  \newenvironment{#1text}{\color{#2}}{\color{black}}
}

\editor{IT}{blue}

\graphicspath{{figs/}}
\usepackage[utf8]{inputenc}
\usepackage[T1]{fontenc}
\usepackage{mathptmx}
\usepackage{etoolbox}

\makeatletter
\def\@email#1#2{%
 \endgroup
 \patchcmd{\titleblock@produce}
  {\frontmatter@RRAPformat}
  {\frontmatter@RRAPformat{\produce@RRAP{*#1\href{mailto:#2}{#2}}}\frontmatter@RRAPformat}
  {}{}
}%
\makeatother
\begin{document}

\title[Randomized Block Davidson Eigensolvers for Plane-Wave DFT]{Randomized Block Davidson Eigensolvers for Plane-Wave Density-Functional Theory}

\author{Moritz Gubler}
\email{moritz.gubler@psi.ch, taejun.park@epfl.ch, iurii.timrov@psi.ch, laura.grigori@epfl.ch}
\affiliation{PSI Center for Scientific Computing,
Theory, and Data, Paul Scherrer Institute, 5232 Villigen PSI, Switzerland}

\author{Taejun Park}
\affiliation{Institute of Mathematics, \'Ecole Polytechnique F\'ed\'erale de Lausanne (EPFL), CH-1015 Lausanne, Switzerland}

\author{Augustin Bussy}
\affiliation{Swiss National Supercomputing Centre (CSCS), ETH Zurich, CH-6900 Lugano, Switzerland}

\author{Johannes Laute}
\affiliation{PSI Center for Scientific Computing,
Theory, and Data, Paul Scherrer Institute, 5232 Villigen PSI, Switzerland}

\author{Nicola Marzari}
\affiliation{Theory and Simulation of Materials (THEOS), and National Centre for Computational Design and Discovery of Novel Materials (MARVEL), \'Ecole Polytechnique F\'ed\'erale de Lausanne (EPFL), CH-1015 Lausanne, Switzerland}
\affiliation{PSI Center for Scientific Computing,
Theory, and Data, Paul Scherrer Institute, 5232 Villigen PSI, Switzerland}

\author{Iurii Timrov}
\affiliation{PSI Center for Scientific Computing,
Theory, and Data, Paul Scherrer Institute, 5232 Villigen PSI, Switzerland}

\author{Laura Grigori}
\affiliation{Institute of Mathematics, \'Ecole Polytechnique F\'ed\'erale de Lausanne (EPFL), CH-1015 Lausanne, Switzerland}
\affiliation{PSI Center for Scientific Computing, Theory, and Data, Paul Scherrer Institute, 5232 Villigen PSI, Switzerland}

\date{\today}

\begin{abstract}
Iterative diagonalization is the dominant cost of plane-wave density-functional theory (DFT), with search-space orthogonalization scaling particularly quickly with problem size and the number of target states. We present a randomized block Davidson-type eigensolver that replaces Euclidean orthogonalization with randomized Gram-Schmidt in a sketched inner product, requiring only a single pass over the basis while keeping its conditioning bounded independently of the input vectors. This modification changes only the Rayleigh-Ritz step, which becomes a definite generalized Hermitian eigenproblem. Ritz extraction remains exact, preserving true Ritz pairs and the interlacing property that makes each band energy an upper bound on the true one.
The method is implemented in mixed precision for CPUs and GPUs from a single Julia code, interfaces matrix-free with DFTK, and is released in the open-source RandESC library. On sparse test problems with a fixed number of eigenpairs, the sketched solver overtakes its deterministic counterpart beyond matrix dimensions of about $2\times 10^4$ and is $25\%$ faster at $5\times 10^5$. In full self-consistent field DFT calculations, the benefits of sketching depend on the application case. With standard pseudopotentials and plane-wave cutoffs, efficiency gains are limited. As the number of requested states grows with system size, orthogonalization savings are offset by the generalized eigenproblem. However, some highly accurate pseudopotentals require larger cutoffs, and in this regime, the randomized solver is measurably better by 6\% overall. Its relative efficiency also improves as simulated systems get larger. 
\end{abstract}

\maketitle

\section{Introduction}
\label{sec:Intro}

Density-functional theory (DFT)~\cite{Hohenberg:1964, KohnSham:1965} predicts the electronic structure of molecules and materials from first principles and, by a wide margin, is the most heavily used method in computational materials science. DFT calculations account for a substantial share of the workload at major supercomputing centres~\cite{Marzari:2021}. In the Kohn-Sham formulation of DFT, the charge density of $\nel$ interacting electrons $\rho(\vb r) = \sum_{i=1}^{\nel} \abs{\psi_i(\vb r)}^2$ is reproduced by $\nel$ auxiliary, non-interacting Kohn-Sham orbitals $\{\psi_i\}$ that minimize the energy functional
\begin{equation}
    E[\{\psi_i\}] = T_s[\{\psi_i\}] + E_{en}[\rho] + E_{ee}[\rho] + E_{xc}[\rho]
    \label{eq:ks_energy}
\end{equation}
subject to the orthonormality constraint $\langle \psi_i | \psi_j \rangle = \delta_{ij}$. Here, $T_s$ is the kinetic energy of the non-interacting system, $E_{en}$ and $E_{ee}$ are the electrostatic energy of the electrons with the nuclei and with each other, respectively, and $E_{xc}$ the exchange-correlation energy. The energy functional in Eq.~\eqref{eq:ks_energy} can be minimized directly, on the manifold of orthonormal orbitals, without ever forming an eigenvalue problem~\cite{Payne:1992}. The route taken here, and the standard one in electronic-structure codes based on plane-wave basis sets, is instead to enforce the constraint through Lagrange multipliers and set the constrained gradient of Eq.~\eqref{eq:ks_energy} to zero. This yields the Kohn-Sham equations
\begin{equation}
    H[\rho]\, \psi_i = \varepsilon_i\, \psi_i, 
    \label{eq:ks_eq}
\end{equation}
where
\begin{equation}
     H[\rho] = -\tfrac{1}{2}\nabla^2 + v_{\mathrm{KS}}[\rho] ,
    \label{eq:ks_Hamiltonian}
\end{equation}
is the Hamiltonian of the electronic system. Here, $\varepsilon_i$ are the eigenvalues and $\psi_i$ are the eigenstates of the Eq.~\eqref{eq:ks_eq}, while $v_{\mathrm{KS}}$ is the Kohn-Sham potential. Since $v_{\mathrm{KS}}[\rho]$ depends on the density it is meant to determine, Eq.~\eqref{eq:ks_eq} is solved self-consistently: at iteration $n$, it is diagonalized at the current density $\rho_n$ for the $\nel$ lowest eigenpairs, a new density is formed by mixing $\sum_i \abs{\psi_i}^2$ into $\rho_n$, and the cycle repeats until the density stops changing with a desired threshold. In a plane-wave basis, $H$ is never assembled explicitly; its action on a vector is evaluated matrix-free via fast Fourier transforms (FFTs) between real and reciprocal space, and diagonalizing it for the occupied orbitals is the dominant cost of every self-consistent field (SCF) step.

The size of Eq.~\eqref{eq:ks_eq} sets the numerical problem this paper addresses. Discretized in a plane-wave basis of size $\npw$, $H$ is an $\npw \times \npw$ Hermitian matrix, of which only the $\nel \ll \npw$ lowest eigenpairs are wanted. $\npw$ grows with the size and the plane-wave kinetic-energy cutoff of the system and routinely reaches into the hundreds of thousands. Diagonalizing $H$ densely, as LAPACK or ScaLAPACK would, costs $O(\npw^3)$ and returns all $\npw$ eigenpairs, the overwhelming majority of which are never used. Both the cost and the memory this requires are prohibitive once $\npw$ moves beyond a few thousand. What is desired instead is a method that exploits the matrix-free structure of $H$ and returns only the eigenpairs asked for. Krylov subspace methods such as the Lanczos algorithm~\cite{Lanczos:1950} build such a space from repeated applications of $H$; preconditioned variants, among them LOBPCG~\cite{Knyazev:2001} and the Davidson family~\cite{Davidson:75, Sleijpen:00} taken up in Sec.~\ref{sec:JD}, additionally use an approximate inverse of $H - \theta \vb I$ to accelerate convergence; Chebyshev-filtered subspace iteration~\cite{Zhou:2006} instead applies a polynomial filter to $H$ to damp the unwanted part of the spectrum before projecting. Every such method builds a low-dimensional space $\mathcal{V}_j$ of dimension $j \ll \npw$ and replaces the $\npw$-dimensional problem of Eq.~\eqref{eq:ks_eq} by its projection onto $\mathcal{V}_j$, a $j\times j$ problem whose lowest eigenpairs approximate those of $H$ with an accuracy set by how well $\mathcal{V}_j$ captures the wanted invariant subspace. For orthonormal $\vb V_j$, this is a $j\times j$ Hermitian eigenproblem whose $k$ lowest eigenpairs, the Ritz pairs, approximate those of $H$, with the accuracy set by how well $\mathcal{V}_j$ captures the wanted invariant subspace. All of them share the same structure: a search space built from matrix-vector products and refined by repeated Rayleigh-Ritz projections, which the rest of this section develops in general with $n = \npw$ and $k = \nel$.

The remainder of this paper is organized as follows. Sections~\ref{sec:JD} and~\ref{sec:Sketching} review the two known ingredients that this work combines, the generalized Davidson algorithm and randomized sketching. Section~\ref{sec:JD_Sketching} presents the new algorithm that results from orthogonalizing the Davidson search space in a sketched inner product, and Sec.~\ref{sec:Implementation} its mixed-precision CPU and GPU implementation in the open source \texttt{RandESC} library~\cite{RandESC_repo}, interfaced with the plane-wave DFT code DFTK~\cite{Herbst:2021}. It is then benchmarked (Sec.~\ref{sec:Benchmarks}) in three regimes: sparse test matrices with a fixed number of wanted eigenpairs, where sketching gives a speed-up that grows with the matrix dimension, and two sets of SCF DFT calculations on crystal structures from the MC3D database~\cite{Huber:26}, where the number of wanted states grows with the system and the benefit of sketching is bounded. In the first set of DFT calculations, standard pseudopotentials and plane-wave cutoffs are employed. In the second set, the use of highly accurate pseudopotentials requiring very high cutoffs is investigated.

\subsection{Generalized Davidson algorithm}
\label{sec:JD}

We consider the Hermitian eigenvalue problem
\begin{equation}
    \vb A \vb x_i = \lambda_i \vb x_i, \qquad \lambda_1 \le \lambda_2 \le \dots \le \lambda_n ,
    \label{eq:evp}
\end{equation}
with $\vb A \in \mathbb{C}^{n \times n}$ known only through its action on a vector $x_i$, and we ask for the $k \ll n$ lowest eigenpairs. Throughout, $\vb X^\dagger$ denotes the conjugate transpose of $\vb X$, $\norm{\cdot}$ the Euclidean norm, and $\kappa(\cdot)$ the spectral condition number.

The eigensolvers considered in this work, including the LOBPCG~\cite{Knyazev:2001} reference of Sec.~\ref{sec:Benchmarks}, are subspace projection methods. They maintain a search space $\mathcal V_j = \mathrm{span}(\vb V_j)$ with $\vb V_j \in \mathbb{C}^{n \times j}$ and $k \le j \ll n$, and extract approximate eigenpairs from it by the Rayleigh-Ritz procedure: with $\vb V_j$ orthonormal one forms the projected matrix $\vb H_j = \vb V_j^\dagger \vb A \vb V_j$, solves the small dense Hermitian eigenproblem
\begin{equation}
    \vb H_j \vb y_i = \theta_i \vb y_i ,
    \label{eq:rr}
\end{equation}
and retains the $k$ lowest Ritz pairs $(\theta_i, \vb u_i = \vb V_j \vb y_i)$ together with their residuals
\begin{equation}
    \vb r_i = \vb A \vb u_i - \theta_i \vb u_i .
    \label{eq:residual}
\end{equation}
Because $\vb A$ is Hermitian the Ritz values interlace the exact eigenvalues, $\lambda_i \le \theta_i$, so that the band energies are approached from above, and decrease monotonically between restarts as $\mathcal V_j$ grows. What distinguishes the individual methods is not this extraction step, which they share, but the rule by which $\mathcal V_j$ is expanded.

Davidson's original method~\cite{Davidson:75} expands by the diagonally preconditioned residual,
\begin{equation}
    \vb t_i = (\vb D - \theta_i \vb I)^{-1} \vb r_i, \qquad \vb D = \mathrm{diag}(\vb A) ,
    \label{eq:davidson}
\end{equation}
and appends $\vb t_i$ to $\mathcal V_j$ after orthogonalization; the block version that treats several roots at once is due to Liu~\cite{Liu:78}. The method was designed for configuration-interaction matrices, which are strongly diagonally dominant, so that $\vb D - \theta_i \vb I$ is a good and essentially free approximation of $\vb A - \theta_i \vb I$. The step that frees the method from the assumption of diagonal dominance was taken by Morgan and Scott~\cite{Morgan:86}, who observed that nothing in the derivation requires the diagonal specifically. Replacing it by any cheaply applicable $\vb K_i$ for $\vb A - \theta_i \vb I$ gives the \emph{generalized Davidson} (GD) expansion
\begin{equation}
    \vb t_i = \vb K_i^{-1} \vb r_i ,
    \label{eq:gd}
\end{equation}
which is the form used in this work. Generalized Davidson is thus preconditioned inverse iteration accelerated by a growing subspace, and it inherits the entire body of preconditioning technique developed for linear systems~\cite{Crouzeix:94, Saad:2011}; near-optimal restarted variants and robust implementations are surveyed in \cite{Stathopoulos:98,Stathopoulos:2010}. LOBPCG~\cite{Knyazev:2001} uses the same expansion vectors \eqref{eq:gd} but restarts at every step onto the fixed $3k$-dimensional space spanned by the current Ritz vectors, the preconditioned residuals and the previous search directions, in place of the accumulated history that generalized Davidson keeps; this is the sense in which the two are compared in Sec.~\ref{sec:Benchmarks}.

One property of \eqref{eq:gd} governs the choice of $\vb K_i$: the expansion degrades as the preconditioner improves. Nominally $\vb K_i$ approximates $\vb A - \theta_i \vb I$, but were the approximation exact then $\vb t_i = (\vb A - \theta_i \vb I)^{-1}(\vb A - \theta_i \vb I) \vb u_i = \vb u_i \in \mathcal V_j$, the expansion vector would already lie in the search space, and the method would stall~\cite{Crouzeix:94}. In practice the expansion vector is orthogonalized against $\mathcal V_j$ before being appended (Sec.~\ref{sec:correction}), which removes the component of $\vb K_i^{-1} \vb r_i$ already in the search space and leaves only its orthogonal part; the projection that Jacobi-Davidson \cite{Sleijpen:00} builds into the correction equation is thus supplied here by the orthogonalization step, and \eqref{eq:gd} is used as it stands.

Two ingredients make this practical. First, the search space cannot grow without bound, since memory $O(nj)$, orthogonalization $O(nj)$ per appended vector and the $O(j^3)$ solve \eqref{eq:rr} all grow with $j$; it is restarted at $j_{\max}$ onto the span of the $j_{\min} < j_{\max}$ best Ritz vectors~\cite{Stathopoulos:98}. Second, the $k$ bands do not converge at the same rate. Once $\norm{\vb r_i}$ falls below the tolerance, band $i$ is dropped from the active set $\mathcal A$ and receives no further correction, so that each iteration applies the Hamiltonian only $\abs{\mathcal A}$ times. The locking is \emph{soft}~\cite{Knyazev:2004,Knyazev:2007} as the converged Ritz vector $\vb u_i$ is left in the search space rather than deflated out of it, and the remaining bands continue to be extracted from a space that contains it.

Finally, we record the cost of one iteration, since this is what randomization acts on. With $k$ active bands and $\dim \mathcal V_j = j$, the solver performs (i) $k$ Hamiltonian applications, each a pair of fast Fourier transforms plus the nonlocal projection, at a cost independent of $j$; (ii) the orthogonalization of the $k$ new vectors against $\vb V_j$, at $O(njk)$; (iii) the update of $\vb H_j$, at $O(njk)$; (iv) the dense solve \eqref{eq:rr}, at $O(j^3)$ and independent of $n$; and (v) the formation of the Ritz vectors $\vb V_j \vb Y$ and their residuals $\vb W_j \vb Y$, at $O(njk)$. Only (ii), (iii) and (v) couple $n$ to $j$, and (ii) may be performed twice for stability. The second pass is needed for a reason specific to Davidson methods: as the residuals shrink, the corrections $\vb K_i^{-1} \vb r_i$ become nearly linearly dependent on the space that produced them, and a single pass of classical Gram-Schmidt loses orthogonality on such a block. The resulting two-pass scheme, classical Gram-Schmidt with reorthogonalization (CGS2), costs $4njk$ per pass and so $8njk$ in all~\cite{Giraud:2005}. Randomization reduces this to $2njk$ in a single pass (Sec.~\ref{sec:Sketching}) while still delivering a well-conditioned basis.


\subsection{Sketching}
\label{sec:Sketching}

Of the operations whose cost grows with the search space, orthogonalization is the one randomization acts on. The Hamiltonian application is not among them, since it is a batch of independent transforms, one per active band, at a cost independent of the current dimension $j$. The orthogonalization, by contrast, forms $\vb V_j^\dagger \vb t$ against the whole basis, and CGS2 pays for it twice; on distributed hardware each pass is also a synchronization point, which is the usual motivation for randomizing it~\cite{Balabanov:2022}. What we want is to reduce that cost without giving up the conditioning the second pass was there to secure, and the tool for it is the subspace embedding.

\emph{Subspace embeddings.} A matrix $\vb S \in \mathbb{C}^{s \times n}$ with $s \ll n$ is a subspace embedding for a subspace $\mathcal U \subset \mathbb{C}^n$ with distortion $\eta \in (0,1)$ if
\begin{equation}
    (1-\eta) \norm{\vb v}^2 \le \norm{\vb S \vb v}^2 \le (1+\eta) \norm{\vb v}^2
    \qquad \text{for all } \vb v \in \mathcal U .
    \label{eq:embedding}
\end{equation}
Such an $\vb S$ compresses $\mathbb{C}^n$ to $\mathbb{C}^s$ while preserving the geometry of $\mathcal U$ to within the factor $1 \pm \eta$. By the polarization identity, \eqref{eq:embedding} also controls inner products,
\begin{equation}
    \abs{\langle \vb S \vb u, \vb S \vb v \rangle - \langle \vb u, \vb v \rangle} \le \eta \norm{\vb u} \norm{\vb v} ,
    \label{eq:sketched_ip}
\end{equation}
so that the \emph{sketched inner product} $\langle \vb u, \vb v \rangle_{\vb S} := \langle \vb S \vb u, \vb S \vb v \rangle$ and its induced norm $\norm{\vb u}_{\vb S} := \norm{\vb S \vb u}$ agree with their Euclidean counterparts on $\mathcal U$ to within $\eta$, at a cost of $O(s)$ rather than $O(n)$ once the sketches are available.

The subspace $\mathcal U$ is usually not known in advance. One therefore draws $\vb S$ at random from a distribution for which \eqref{eq:embedding} holds, with probability at least $1 - \delta$, for \emph{any} fixed subspace of a given dimension $d$; such an $\vb S$ is an oblivious subspace embedding~\cite{Sarlos:06, Woodruff:2014, Martinsson:2020}. The classical construction is a Gaussian matrix, for which $s = O(\eta^{-2}(d + \log(1/\delta)))$ suffices~\cite{Halko:2011}, but applying a dense $s \times n$ matrix costs $O(sn)$ per vector and would defeat the purpose. Two cheaper families are used in practice: subsampled randomized trigonometric transforms~\cite{Ailon:2009, Tropp:2011}, which apply a fast transform and then select $s$ rows at random, at cost $O(n \log n)$; and sparse sign embeddings~\cite{Clarkson:2017, Nelson:2013, Cohen:2016}, whose columns carry only $\zeta$ nonzero entries $\pm 1/\sqrt{\zeta}$ in random positions, at cost $O(\zeta n)$. The recommended value of $\zeta$ in practice is a small constant, say $4$ or $8$.

\emph{Sketch-orthonormal bases.} The use we make of \eqref{eq:embedding} is not to approximate the projected eigenproblem but to condition the basis. A matrix $\tilde{\vb V} \in \mathbb{C}^{n \times j}$ is called sketch-orthonormal if
\begin{equation}
    (\vb S \tilde{\vb V})^\dagger (\vb S \tilde{\vb V}) = \vb I_j ,
    \label{eq:sketch_orth}
\end{equation}
that is, its sketch $\vb Q := \vb S \tilde{\vb V} \in \mathbb{C}^{s \times j}$ is orthonormal. If $\vb S$ is an $\eta$-embedding for $\mathrm{span}(\tilde{\vb V})$, then applying \eqref{eq:embedding} to $\vb v = \tilde{\vb V} \vb y$ gives $(1-\eta) \norm{\tilde{\vb V} \vb y}^2 \le \norm{\vb y}^2 \le (1+\eta) \norm{\tilde{\vb V} \vb y}^2$ for every $\vb y$, so that all singular values of $\tilde{\vb V}$ lie in $[(1+\eta)^{-1/2}, (1-\eta)^{-1/2}]$ and
\begin{equation}
    \kappa(\tilde{\vb V}) \le \sqrt{\frac{1+\eta}{1-\eta}} .
    \label{eq:cond}
\end{equation}
A sketch-orthonormal basis is not orthonormal, but at $\eta = 1/2$ its condition number is at most $\sqrt 3$, regardless of how nearly linearly dependent the input vectors are. That is the key point. Orthonormality was never the goal in itself; it was a means of obtaining a well-conditioned basis, and \eqref{eq:sketch_orth} reaches it at a fraction of the cost.

\emph{Randomized orthogonalization.} The basis is built by the randomized orthogonalization process of Balabanov and Grigori~\cite{Balabanov:2022, Balabanov:2021}, in the randomized classical Gram-Schmidt form (RCGS) of Algorithm~\ref{alg:rgs}, which proceeds in two stages. The block $\vb T$ of new vectors is sketched once, $\vb P = \vb S \vb T$; its coefficients against the existing basis are formed in the sketched space, $\vb H = \vb Q^\dagger \vb P$, and removed from both spaces at once,
\begin{equation}
    \hat{\vb T} = \vb T - \tilde{\vb V} \vb H , \qquad \hat{\vb P} = \vb P - \vb Q \vb H ,
    \label{eq:rgs}
\end{equation}
so that $\hat{\vb P}$ is the sketch of $\hat{\vb T}$ and need not be recomputed. The block is then orthonormalized against itself in the sketched norm. The first stage, which carries the dominant cost, is a pair of matrix-matrix products; the second, of cost $O(nk^2)$, is a column loop in the default variant and a single batched factorization in the sketched-QR and sketched-Cholesky-QR variants, which are preferred on GPUs. In the latter the Cholesky factor is taken of $\hat{\vb P}^\dagger \hat{\vb P}$ rather than $\hat{\vb T}^\dagger \hat{\vb T}$, so even that factorization costs $O(sk^2)$ instead of $O(nk^2)$.

Orthogonalization and extraction are independent places to randomize. The process above sketches only the construction of the basis, and leaves the Rayleigh-Ritz step of Sec.~\ref{sec:JD} to be taken in the full space. The alternative, a randomized Rayleigh-Ritz in which the projected matrices are themselves formed from sketches, is cheaper again, as in \cite{Nakatsukasa:2024}; for Hermitian $\vb A$ it costs more than it saves, and Sec.~\ref{sec:extraction} takes the other route.

Two caveats attach to \eqref{eq:embedding} inside an iterative method. The search space depends on $\vb S$, so the oblivious guarantee, which holds for a subspace fixed in advance, does not apply verbatim; as in randomized Krylov solvers~\cite{Balabanov:2022, Nakatsukasa:2024}, we assume $\vb S$ embeds the spaces actually generated. Because the Rayleigh-Ritz step is taken in the full space, $\vb A \vb{\widetilde{V}}_j$ is never sketched, and $\vb S$ need embed no more than $\operatorname{span}([\vb{\widetilde{V}}_j,\vb T])$, of dimension at most $j_{\max}$ once the block is appended. A single draw with $s$ tied to $j_{\max}$ suffices at every $j \le j_{\max}$, and no resketching is needed. The assumption is at least checkable at a cost independent of $n$: with $\vb G_j = \tilde{\vb V}_j^\dagger \tilde{\vb V}_j$, \eqref{eq:embedding} holds on $\mathrm{span}(\tilde{\vb V}_j)$ if and only if the eigenvalues of $\vb G_j $ lie in $[(1+\eta)^{-1},\, (1-\eta)^{-1}]$, and the algorithm of Sec.~\ref{sec:JD_Sketching} forms $\vb G_j$ in any case.

\begin{algorithm}[H]
\caption{$\textsc{rcgs}$: sketched orthogonalization of a block $\vb T$ against $\tilde{\vb V}$}
\label{alg:rgs}
\begin{algorithmic}[1]
\Require $\vb T \in \mathbb{C}^{n \times k_{\mathrm a}}$; sketch $\vb S$; basis $\tilde{\vb V} \in \mathbb{C}^{n \times j}$ with $\vb Q = \vb S \tilde{\vb V}$
\Ensure $\hat{\vb T}$ and its sketch $\hat{\vb P} = \vb S \hat{\vb T}$, sketch-orthonormal and sketch-orthogonal to $\tilde{\vb V}$
\State $\vb P \gets \vb S \vb T$ \Comment{$\zeta n k_{\mathrm a}$, independent of $j$}
\For{$p = 1, \dots, p_{\mathrm s}$} \Comment{stage 1, blocked; $p_{\mathrm s} = 1$ (\texttt{rcgs}) or $2$ (\texttt{rcgs2})}
  \State $\vb H \gets \vb Q^\dagger \vb P$ \Comment{$2sjk_{\mathrm a}$, independent of $n$}
  \State $\vb T \gets \vb T - \tilde{\vb V} \vb H$, \quad $\vb P \gets \vb P - \vb Q \vb H$ \Comment{$2njk_{\mathrm a}$: the only $n$--$j$ term}
\EndFor
\For{$\ell = 1, \dots, k_{\mathrm a}$} \Comment{stage 2, column loop; batched in the \texttt{rqr}/\texttt{cholqr} variants}
  \State $\vb h \gets \hat{\vb P}_{1:\ell-1}^\dagger \vb p_\ell$ 
  \State $\vb t_\ell \gets \vb t_\ell - \hat{\vb T}_{1:\ell-1} \vb h$, \quad $\vb p_\ell \gets \vb p_\ell - \hat{\vb P}_{1:\ell-1} \vb h$
\EndFor
\end{algorithmic}
\end{algorithm}

\section{Sketching the generalized Davidson algorithm}
\label{sec:JD_Sketching}

The randomized solver follows from Sec.~\ref{sec:JD} by a single substitution: the basis is orthogonalized in the sketched inner product rather than the Euclidean one, by Algorithm~\ref{alg:rgs}. The rest of this section works out what that entails. The method is summarized in Algorithm~\ref{alg:rgd}.
 
\subsection{Extraction}
\label{sec:extraction}
 
The Rayleigh-Ritz step \eqref{eq:rr} used the orthonormality of $\vb V_j$ twice: to write the Galerkin condition $\vb V_j^\dagger (\vb A - \theta \vb I) \vb V_j \vb y = 0$ as a standard eigenproblem, and to make $\vb u_i = \vb V_j \vb y_i$ normalized. A sketch-orthonormal basis satisfies neither, and how one responds is the central design decision of the method.
 
The first option, raised in Sec.~\ref{sec:Sketching}, is to sketch the projection too. With $\vb Q_j = \vb S \tilde{\vb V}_j$ and $\vb B_j = \vb S \vb A \tilde{\vb V}_j$, sketch-orthonormality \eqref{eq:sketch_orth} turns the Galerkin condition into the standard eigenproblem $\vb Q_j^\dagger \vb B_j \tilde{\vb y} = \tilde\theta \tilde{\vb y}$, which costs $O(sjk)$ per iteration to update and needs no inner products against the full-space basis at all. However, for Hermitian $\vb A$, it gives up more than it saves: $\vb Q_j^\dagger \vb B_j$ is not Hermitian, so its eigenvalues are complex in general, its Ritz vectors are not orthonormal, and solving a non-Hermitian eigenvalue problem can be far more expensive than a Hermitian one \cite{Golub:2013}. For non-Hermitian $\vb A$ this concern does not apply, since the projected matrix is not Hermitian in any case, and sketching the projection is then the natural choice; a randomized Jacobi-Davidson method of this kind is analyzed in \cite{Grigori:2026}. In either case, one cost remains. Sketching $\vb A \tilde{\vb V}_j$ alongside $\tilde{\vb V}_j$ means $\vb S$ must embed $\operatorname{span}([\tilde{\vb V}_j, \vb A \tilde{\vb V}_j])$, of up to twice the dimension, so the $O(sjk)$ above is measured against a larger $s$.
 
The second option, taken here, keeps the projection exact and pays instead for the non-orthonormal basis. The Galerkin condition on $\mathrm{span}(\tilde{\vb V}_j)$ reads
\begin{equation}
    \vb H_j \vb y_i = \theta_i \vb G_j \vb y_i, \qquad
    \vb H_j = \tilde{\vb V}_j^\dagger \vb A \tilde{\vb V}_j, \quad
    \vb G_j = \tilde{\vb V}_j^\dagger \tilde{\vb V}_j ,
    \label{eq:ghep}
\end{equation}
a generalized Hermitian eigenvalue problem in place of the standard one \eqref{eq:rr}. Both matrices are Hermitian, $\vb G_j$ is positive definite, and \eqref{eq:cond} bounds $\kappa(\vb G_j) \le (1+\eta)/(1-\eta)$, so the pencil is definite and reduces to standard form by a Cholesky factorization of $\vb G_j$. The extraction is then exactly that of Sec.~\ref{sec:JD}, the $\theta_i$ are true Rayleigh quotients of true Ritz vectors, they interlace, and the classical Davidson theory carries over unchanged. Randomization is confined to the construction of the basis and does not touch the approximation itself; all that survives of it is the departure of $\vb G_j$ from the identity, which \eqref{eq:cond} bounds.
 
One consequence deserves separate mention. Since $\vb Y^\dagger \vb G_j \vb Y = \vb I$, the Ritz vectors $\vb U = \tilde{\vb V}_j \vb Y$ satisfy $\vb U^\dagger \vb U = \vb I$: they are exactly orthonormal though the basis they were formed from is not, so nothing is owed at the end.

\subsection{Correction and expansion}
\label{sec:correction}
 
The corrections are those of Sec.~\ref{sec:JD} and are untouched by randomization: the preconditioner is applied to the block of active residuals, $\vb T = \vb K^{-1} \vb R_{\mathcal A}$, and the result passed to Algorithm~\ref{alg:rgs}. Randomization pays only where many inner products are taken against a common basis, and that is the orthogonalization.

\subsection{Restart}
 
The exact orthonormality of the Ritz vectors makes the restart nearly free. Retaining the span of the $j_{\min}$ lowest means rotating by $\vb Y_{1:j_{\min}}$, and since $\vb Y^\dagger \vb G_j \vb Y = \vb I$ the rotated basis is orthonormal, so $\vb H_{j_{\min}} = \mathrm{diag}(\theta_1, \dots, \theta_{j_{\min}})$ and $\vb G_{j_{\min}} = \vb I$ are known without recomputation. The rotations of $\tilde{\vb V}_j$, $\vb W_j$ and $\vb Q_j$ are those the deterministic solver performs anyway, no resketching is needed, and the retained Ritz vectors stay in the span, so soft locking is unaffected. The restarted basis is orthonormal and therefore no longer sketch-orthonormal, which the method tolerates because \eqref{eq:ghep} never assumed that it was; the one consequence is that the sketched projection \eqref{eq:rgs} of the next block is oblique rather than orthogonal, by a relative amount bounded by $\eta$, until sketch-orthonormality is re-established.
 
\begin{algorithm}[H]
\caption{Generalized Davidson with sketched orthogonality}
\label{alg:rgd}
\begin{algorithmic}[1]
\Require Hamiltonian action $\vb v \mapsto \vb A \vb v$; preconditioner $\vb K^{-1}$; initial block $\vb V_0 \in \mathbb{C}^{n \times k}$; tolerance $\tau$; restart dimensions $j_{\min} < j_{\max}$; sketch size $s$
\Ensure $k$ Ritz pairs $(\theta_i, \vb u_i)$ with $\vb U^\dagger \vb U = \vb I$ and $\norm{\vb A \vb u_i - \theta_i \vb u_i} \le \tau$
\State Draw $\vb S \in \mathbb{C}^{s \times n}$ \Comment{sparse sign, $\zeta = 4$, $s = 5 j_{\max}$}
\State $(\tilde{\vb V}, \vb Q) \gets \textsc{rcgs}(\vb V_0, \vb S)$; \; $\vb W \gets \vb A \tilde{\vb V}$; \; $j \gets k$
\Repeat
  \State $\vb G \gets \tilde{\vb V}^\dagger \tilde{\vb V}$, \quad $\vb H \gets \tilde{\vb V}^\dagger \vb W$
  \State Solve $\vb H \vb Y = \vb G \vb Y \vb{\Theta}$, keep the $k$ lowest, $\vb Y^\dagger \vb G \vb Y = \vb I$ \Comment{definite GHEP, $\kappa(\vb G) \le \frac{1+\eta}{1-\eta}$}
  \State $\vb U \gets \tilde{\vb V} \vb Y$, \quad $\vb R \gets \vb W \vb Y - \vb U \vb{\Theta}$ \Comment{$\vb U^\dagger \vb U = \vb I$}
  \State $\mathcal A \gets \{ i : \norm{\vb r_i} > \tau \}$ \Comment{soft locking}
  \If{$\mathcal A = \emptyset$}
    \State \Return $(\vb{\Theta}, \vb U)$
  \EndIf
  \State $\vb T \gets \vb K^{-1} \vb R_{\mathcal A}$ \Comment{preconditioned residuals \eqref{eq:gd}, blockwise}
  \State $(\hat{\vb T}, \hat{\vb P}) \gets \textsc{rcgs}(\vb T, \vb S, \tilde{\vb V}, \vb Q)$ \Comment{Algorithm~\ref{alg:rgs}: sketch-orthogonalize $\vb T$ against $\tilde{\vb V}$}
  \State $\tilde{\vb V} \gets [\tilde{\vb V}, \hat{\vb T}]$, \; $\vb Q \gets [\vb Q, \hat{\vb P}]$, \; $\vb W \gets [\vb W, \vb A \hat{\vb T}]$, \; $j \gets j + \abs{\mathcal A}$ 
  \If{$j \ge j_{\max}$}
    \State $\tilde{\vb V} \gets \tilde{\vb V} \vb Y_{1:j_{\min}}$, $\vb W \gets \vb W \vb Y_{1:j_{\min}}$, $\vb Q \gets \vb Q \vb Y_{1:j_{\min}}$, $j \gets j_{\min}$
    \State $\vb H \gets \mathrm{diag}(\theta_1, \dots, \theta_{j_{\min}})$, \; $\vb G \gets \vb I$ \Comment{$\vb Y^\dagger \vb G \vb Y = \vb I$}
  \EndIf
\Until{converged}
\end{algorithmic}
\end{algorithm}

\subsection{Cost}
\label{sec:cost}
 
Randomization changes one term of the iteration and adds one. The deterministic orthogonalization projects the new block out of $\mathcal V_j$ in two passes of $4njk$ (if CGS2 is used); the randomized one does it in a single pass of $2njk$, because the coefficients $\vb H$ are formed in the sketched space and only the update $\tilde{\vb V}_j \vb H$ remains in the full space. Against that saving of $6njk$ the randomized solver forms $\vb G_j$, at $2njk$. The rest is common: $k$ Hamiltonian applications, the update of $\vb H_j$ at $2njk$, the Ritz vectors and residuals at $4njk$, and a dense solve at $O(j^3)$ for which \eqref{eq:ghep} is more expensive than \eqref{eq:rr} but no more expensive in $n$. The totals coupling $n$ to $j$ are $14njk$ (deterministic) against $10njk$ (randomized), and the balance between an $njk$ saving and a $j^3$ price is governed by $n j k \gtrsim j^3$.
 
The ratio $14:10$ may not mean much, because on current hardware neither the saving nor the price is well measured by arithmetic. The randomized pass reads the $n \times j$ basis once, in the single product $\tilde{\vb V}_j \vb H$; the deterministic double pass reads it four times. That array can be large and may not sit in cache, so the factor of four is a factor in memory traffic too, and holds whichever of the two binds. What the sketch touches is smaller still, $s \times j$ or $s \times k$, and single precision serves for it (Sec.~\ref{sec:cpu_gpu_precision}).

The balance therefore turns on how $k$ scales. With $k$ fixed and $n$ growing, the saving grows while the price does not, and a crossover must exist. With $k$ growing in proportion to the problem, as it does when the number of wanted states is set by the system, $j \propto k \propto n$: the two grow at the same rate, and the ratio is fixed by constants that only measurement can supply. The limit in that regime is the $O(j^3)$ dense eigenproblem, which grows as fast as the saving.

\section{Implementation}
\label{sec:Implementation}
\subsection{CPU, GPU, single vs double precision}
\label{sec:cpu_gpu_precision}

The sketched generalized Davidson algorithm of Sec.~\ref{sec:JD_Sketching} is designed to run on both multi-core CPUs and GPUs from a single implementation. As in any plane-wave electronic-structure code, the Hamiltonian is never formed explicitly; its action $\vb A * \vb v$ is evaluated matrix-free through the FFTs between real and reciprocal space. All remaining work, namely maintaining the search space $\vb V$, applying the sketch $\vb S * \vb V$, the projected (Rayleigh-Ritz) eigenproblem, the residuals and the correction step, is expressed in terms of dense block operations acting on the whole block of vectors at once, rather than vector by vector. This block structure is what makes the method well suited to GPUs: the batched FFTs and dense matrix-matrix products expose enough parallelism and arithmetic intensity to use the hardware efficiently, while the same code path runs on the CPU when no accelerator is available. The method is implemented in the Julia programming language, which lets the identical generic code execute on both CPUs and GPUs and carry either single- or double-precision arrays; the implementation is openly available at on GitHub~\cite{RandESC_repo}.

The algorithm operates in mixed precision. We keep the full, $n$-dimensional space in double precision: the search vectors $\vb V$, the Hamiltonian application $\vb A * \vb V$, the residuals and the orthogonalization against the full space are all computed and stored in double-precision complex arithmetic (\texttt{ComplexF64}). The sketched space, by contrast, is handled entirely in single precision: the sketching operator maps double-precision full-space arrays to single-precision sketched arrays, $\vb S:\texttt{ComplexF64}\to\texttt{ComplexF32}$, and every quantity that lives only in the sketched space ($\vb S * \vb V$, $\vb S * \vb A * \vb V$, the small dense projected eigenproblem and the sketched orthogonality constraints) is formed and solved in single precision.

The reason this is admissible is that the sketch is itself already an approximation. The random embedding $\vb S$ only preserves norms and inner products on $\mathrm{span}(\vb V)$ up to a distortion $\eta$ (Sec.~\ref{sec:Sketching}), and this distortion is orders of magnitude larger than the $\sim 10^{-7}$ relative error introduced by single-precision rounding. Reducing the sketched quantities to single precision therefore adds essentially no error beyond what the sketch already incurs, while halving their memory footprint and bandwidth and letting the recurring small dense linear algebra run in single precision. Importantly, the accuracy of the final eigenpairs is not set by the sketch but by the double-precision full space: the residuals $\vb r_i = \vb A * \vb u_i - \theta_i \vb u_i$ used in the convergence test are evaluated in double precision, so the converged eigenpairs reach full double-precision accuracy. The dominant speed-up of the method comes from the dimensionality reduction of sketching itself; carrying the sketched space in single precision makes that reduction essentially free in terms of accuracy.

\subsection{Software}
\label{sec:Software}

All eigensolvers compared in this work are implemented in the open-source \texttt{RandESC} package introduced in Sec.~\ref{sec:cpu_gpu_precision}. The solvers are written against a minimal matrix-free interface: the only problem-specific operations they require are the application of the operator to a block of vectors, $\vb A * \vb V$, and the application of a preconditioner, $\vb K^{-1} * \vb R$. All remaining operations are dense block linear algebra and are independent of the problem at hand. Because of this design, the same solver code runs unchanged on the sparse test matrices of Sec.~\ref{sec:sparse_benchmark} and on the plane-wave Kohn-Sham Hamiltonians of Sec.~\ref{sec:dft_benchmark}. The package provides the deterministic generalized Davidson algorithm of Sec.~\ref{sec:JD} and its sketched variant of Sec.~\ref{sec:JD_Sketching}. Both operate on blocks of vectors, solve the correction equation using a preconditioner, and restart when the search space becomes too large. Converged eigenpairs are soft-locked: they remain in the search space, but no corrections are computed for them anymore. By default, the sketching operator $\vb S$ is a sparse sign embedding with $\zeta = 4$ nonzero entries per column; Gaussian sketches and subsampled random trigonometric transforms are also implemented. The sketch dimension is set to $s = 5 j_{\text{max}}$ by default, where $j_{\text{max}}$ is the largest dimension the search space can reach before a restart is triggered.

The electronic-structure benchmarks are carried out with DFTK~\cite{Herbst:2021}, a plane-wave DFT code written in Julia. DFTK applies the Hamiltonian matrix-free using the FFTs, runs on CPUs and GPUs, and uses LOBPCG~\cite{Knyazev:2001} as its default eigensolver. The solvers of \texttt{RandESC} act as drop-in replacements for LOBPCG in the SCF loop of DFTK. This allows a direct comparison of LOBPCG, generalized Davidson, and sketched generalized Davidson in which the Hamiltonian application, the preconditioner, and all remaining parts of the SCF iteration are identical. All three solvers use the Teter-Payne-Allan kinetic-energy preconditioner~\cite{Teter:89} provided by DFTK.

\section{Benchmarks}
\label{sec:Benchmarks}

We benchmark the sketched generalized Davidson algorithm in three settings that probe different regimes. In the first setting, a fixed number of $k=500$ eigenpairs of a sparse matrix is computed while the dimension $n$ of the matrix grows. In this regime, the cost of keeping the search space orthogonal grows faster than the cost of the Rayleigh-Ritz step, and sketching is expected to pay off. In the second setting, the solvers are compared inside full SCF calculations with DFTK version v0.7.25. There, the number of requested eigenpairs is proportional to the number of electrons and grows together with the dimension of the Hamiltonian, which changes the balance between the individual components of the solvers. The third setting covers full SCF calculations with DFTK again, but in a regime where the plane-wave cutoff is required to be very high for accuracy, thus leading to $k \ll n$. Namely, we use the MPC\cite{Gulans:2026} pseudopotentials to reach near all-electron accuracy.

\subsection{Sparse matrix benchmark}
\label{sec:sparse_benchmark}

\begin{figure*}
    \centering
    \includegraphics[width=.85\linewidth, page=1]{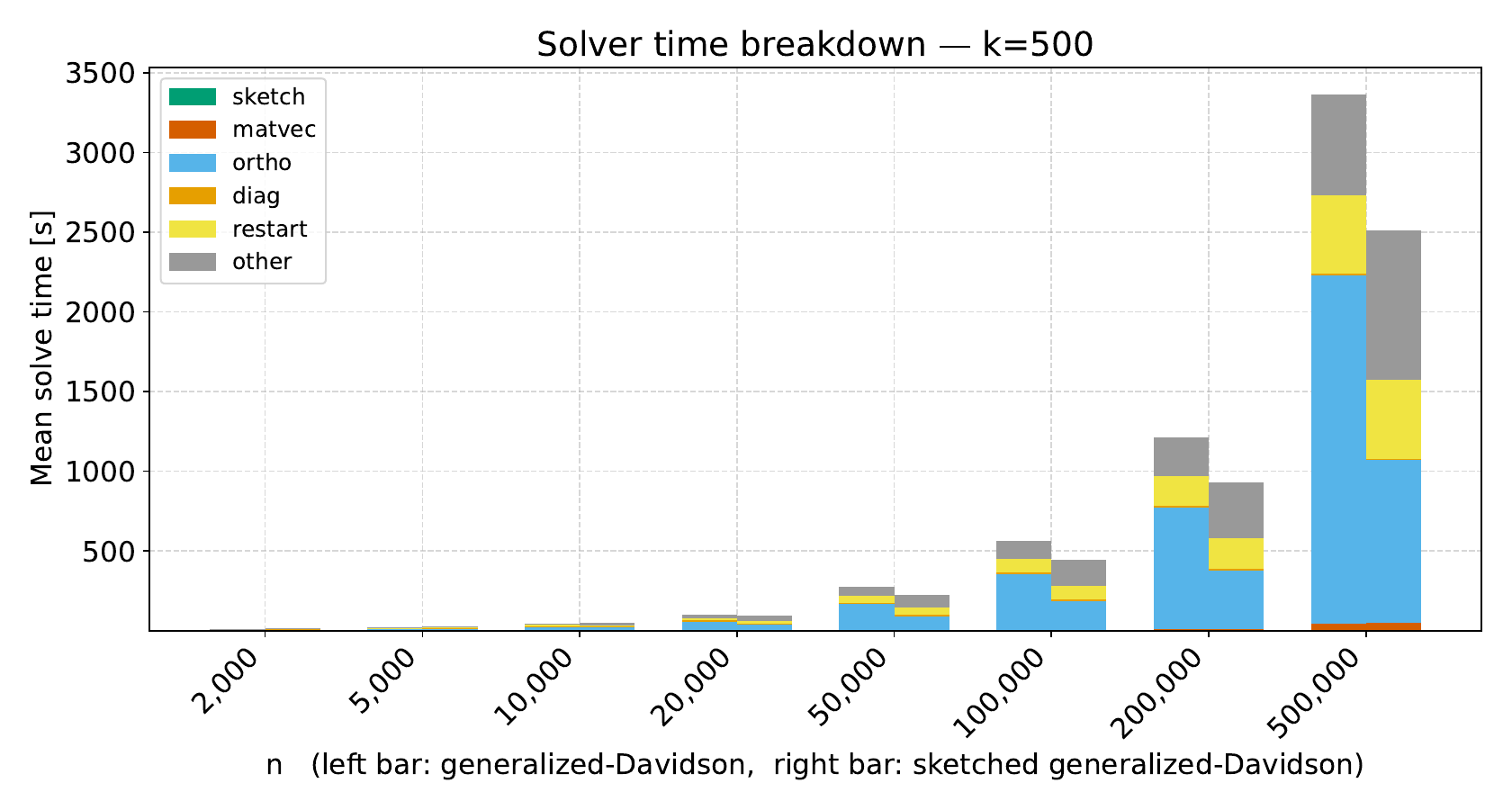}\\
    \includegraphics[width=.85\linewidth, page=2]{figs/sparse_breakdown_from_json.pdf}
    \caption{Breakdown of the solver run time for the sparse matrix benchmark. For each matrix dimension $n$, the left bar corresponds to the deterministic generalized Davidson method and the right bar to its sketched variant. Top: absolute times. Bottom: times relative to the total run time of the deterministic solver at the same $n$.
    Each bar is decomposed into: sketch (sketched variant only), computation of the random sketch $\Theta V$ of the search-space vectors; matvec, the sparse matrix-vector product with the operator $A$; ortho, orthogonalization of the newly generated search-space vectors against the existing basis (standard Gram-Schmidt/QR for the deterministic solver, sketched CGS for the randomized variant); diag, eigendecomposition of the small projected matrix to obtain Ritz pairs; restart, restart step that truncates the search subspace once it exceeds the maximum allowed size; and other, all remaining, individually minor stages.
    }
    \label{fig:sparse_breakdown}
\end{figure*}

In the first benchmark, the lowest $k=500$ eigenpairs of sparse symmetric matrices with dimensions between $n = 2\,000$ and $n = 500\,000$ are computed with the deterministic generalized Davidson method and with its sketched variant. The test matrices are symmetric tridiagonal with a logarithmically growing diagonal,
\begin{equation}
    A_{ii} = \ln(99 + i), \qquad i = 1, \ldots, n,
\end{equation}
and off-diagonal entries drawn uniformly from $(-0.2,\, 0.2)$. The diagonal dominates the off-diagonal entries, so the matrices are positive definite, and the logarithmic ramp gives a slowly growing, bounded spectrum whose condition number remains moderate for all $n$. Both solvers are given the same random initial guesses and the same Jacobi (diagonal) preconditioner, and iterate until the residual norms fall below $10^{-6}$. The reported times are averages over three runs, and just-in-time compilation is excluded from the measurement by a warm-up solve. The benchmark was run on a machine with an Intel Core i9-10940X CPU, where four threads were used to accelerate the linear algebra operations and the remaining work was done serially. The timing breakdowns are shown in Fig.~\ref{fig:sparse_breakdown}.

For small matrices, the sketched solver is slower than the deterministic one. The total run times in Fig.~\ref{fig:sparse_breakdown} cross at $n \approx 2\cdot 10^4$, and from there on the advantage of the sketched solver grows with the matrix dimension. At $n = 5\cdot 10^5$, the sketched solver is roughly $25\%$ faster.

The timing breakdown in Fig.~\ref{fig:sparse_breakdown} explains this behavior. Since the test matrices are very sparse, matrix-vector products are cheap, and the run time of the deterministic solver is increasingly dominated by the orthogonalization of the search space as $n$ grows. The sketched solver replaces this full-space orthogonalization by the sketch-orthogonalization described in Sec.~\ref{sec:JD_Sketching}, which reduces its cost by more than a factor of two at large $n$. The price is a more expensive Rayleigh-Ritz step: because the search space is only sketch-orthonormal, a generalized eigenvalue problem has to be solved instead of a standard one. For small matrices, this additional cost outweighs the savings in the orthogonalization; for large matrices, the orthogonalization savings dominate.

\subsection{Plane-wave DFT benchmarks}
\label{sec:dft_benchmark}

\begin{figure*}
    \centering
    \includegraphics[width=.49\linewidth]{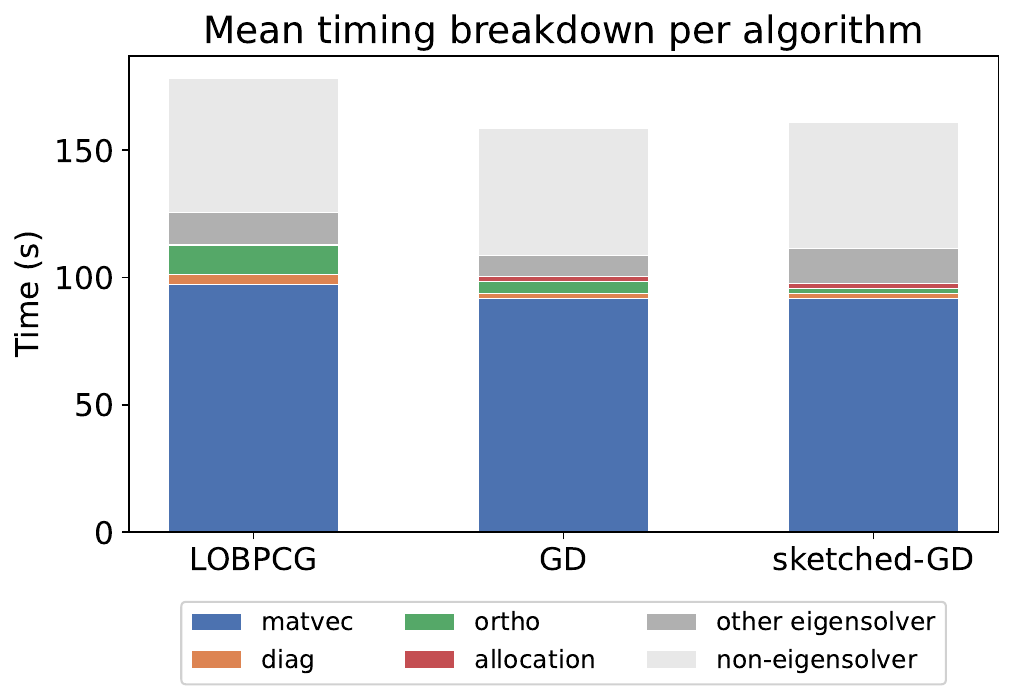}
    \includegraphics[width=.49\linewidth]{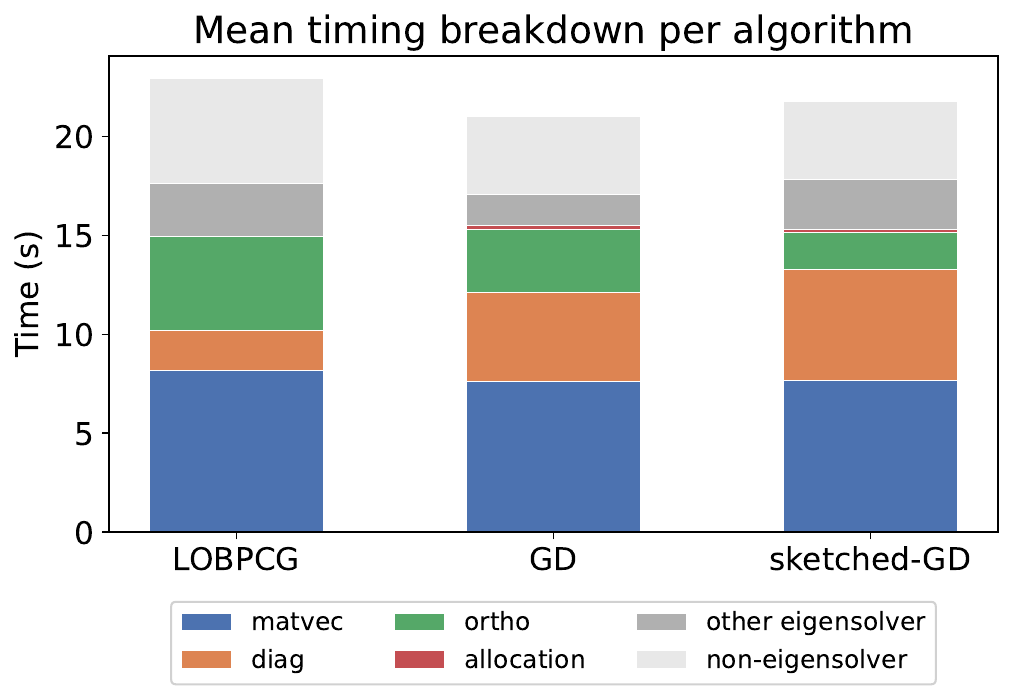}
    \caption{Mean run time per SCF DFT calculation for the three eigensolvers, split into the individual solver components and the solver-independent rest of the SCF iteration. Left: CPU benchmark set. Right: GPU benchmark set.}
    \label{fig:dft_bars}
\end{figure*}

To assess the solvers in a realistic application, complete SCF DFT calculations were performed with DFTK, once with each of the three eigensolvers and otherwise identical settings, for two sets of crystal structures randomly selected from the MC3D database~\cite{Huber:26}:
\begin{itemize}
    \item \textbf{CPU:} 50 non-magnetic structures with 6 to 30 atoms in the unit cell, computed at the $\Gamma$ point with a plane-wave cutoff of 80\,Ry.
    \item \textbf{GPU:} 100 non-magnetic structures with 14 to 64 atoms in the unit cell, computed with a $2\times 2\times 2$ $k$-point grid and a plane-wave cutoff of 90\,Ry.
\end{itemize}

Goedecker-Teter-Hutter~\cite{Goedecker96} Pseudopotentials were used with the PBE~\cite{Perdew:96} exchange-correlation functional.
All calculations were run on the Alps (Daint, GH200 superchip) system at CSCS.

Figure~\ref{fig:dft_bars} shows the mean run time per algorithm, split into the individual solver components and the remaining, solver-independent parts of the SCF iteration. Both generalized Davidson variants are faster than LOBPCG on average: on the CPU, the mean SCF time drops from 178.1\,s with LOBPCG to 158.6\,s with generalized Davidson and 160.8\,s with the sketched variant, and on the GPU from 22.9\,s to 21.0\,s and 21.8\,s, respectively. In all cases, the Hamiltonian application is the largest contribution and is of nearly the same size for the three solvers. The gain of the generalized Davidson variants over LOBPCG comes mostly from the orthogonalization, whose mean cost drops from 11.5\,s with LOBPCG to 4.9\,s and 2.1\,s on the CPU and from 4.7\,s to 3.2\,s and 1.9\,s on the GPU, and on the CPU also from a cheaper Rayleigh-Ritz step. The lower orthogonalization cost of generalized Davidson compared with LOBPCG may appear surprising given its larger search space. However, the DFTK LOBPCG implementation~\cite{Herbst:2021} performs several orthogonality-related operations per iteration: the residual block is projected against the current eigenvectors and search directions, then reorthonormalized with a stability-oriented Cholesky-QR procedure, and coefficient-space orthogonalization is included in the same timer. By contrast, Davidson projects each new correction block against its accumulated basis only once using blocked matrix-matrix products, before orthonormalizing the new block itself. Thus, although Davidson uses a larger subspace, it accesses it fewer times and more efficiently, leading to a lower measured orthogonalization cost.

The deterministic and the sketched generalized Davidson variants perform almost identically in these benchmarks, with the deterministic variant slightly ahead on average. On the one hand, sketching roughly halves the cost of the orthogonalization, as in the sparse benchmark. On the other hand, the sketched solver has to solve a generalized eigenvalue problem in the Rayleigh-Ritz step, which raises its cost from 4.5\,s to 5.6\,s on the GPU, and the application of the sketching operator itself adds overhead. In contrast to the sparse benchmark, where the number of eigenpairs was fixed, in DFT it is proportional to the number of electrons and grows with the system size. The cost of the Rayleigh-Ritz step therefore grows at the same rate as the savings in the orthogonalization, and the two effects nearly cancel.

In the setting presented in this benchmark, both generalized Davidson variants outperform the LOBPCG reference implementation of DFTK, while the benefit of sketching within the generalized Davidson framework is limited by the growth of the Rayleigh-Ritz step with system size.

\subsection{Near all-electron accuracy plane-wave DFT} \label{sec:mpc_pseudos}

For this benchmark, we use the MPC\cite{Gulans:2026} pseudopotentials to reach near all-electron accuracy in the SCF calculations. These pseupotentials are quite hard, and a high plane-wave cutoff is required to resolve the explicit shallow core states. We use the value of 750 Ry, the cutoff recommended by the authors. 50 crystal structures were randomly taken from the MC3D~\cite{Huber:26}, with the following criteria: 24 atoms or more, clean convergence of the SCF cycle at the $\Gamma$-point, and everything fits in the memory of a single GPU of a GH200 superchip on the Daint cluster at CSCS. With such high cutoffs and system sizes, GPU memory becomes the main limitation. DFTK is mostly GPU resident, and at the time of writing, the software does not allow MPI parallelization over the plane waves (i.e. reciprocal-lattice vectors).

\begin{figure}
    \includegraphics[width=.98
    \linewidth]{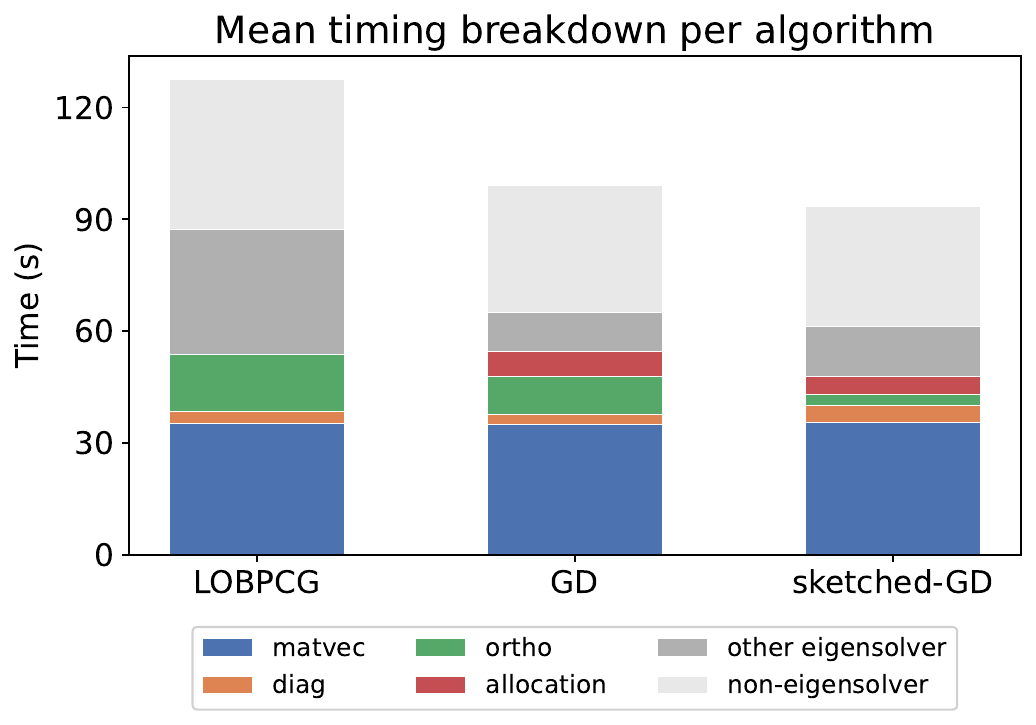}
    \caption{Mean run time per SCF DFT calculation for the three eigensolvers on GPUs, using the MPC pseudopotentials. Timings are split into the individual solver components and the solver-independent rest of the SCF iteration.}
    \label{fig:mpc_bars}
\end{figure}

Figure \ref{fig:mpc_bars} shows the average SCF GPU runtimes over the 50 benchmark systems, for the three solvers of interest. Both flavors of generalized Davidson now clearly outperform the native LOBPCG implementation of DFTK v0.7.25. Additionally, we observe a measurable efficiency advantage of the sketched generalized Davidson solver against its generalized counterpart. When switching sketching on, the solver is 3.7 seconds faster on average. Crucially, the cost of orthogonalization goes down from 10.2 s to 2.8 s, while diagonalization only increases from 2.6 s to 4.7 s. Overall, this represents a 6\% efficiency gain. We also note that in this regime, the eigensolver cost is lower relative to the whole SCF calculation compared to figure \ref{fig:dft_bars}.

\begin{figure*}
    \centering
    \includegraphics[width=.98\linewidth]{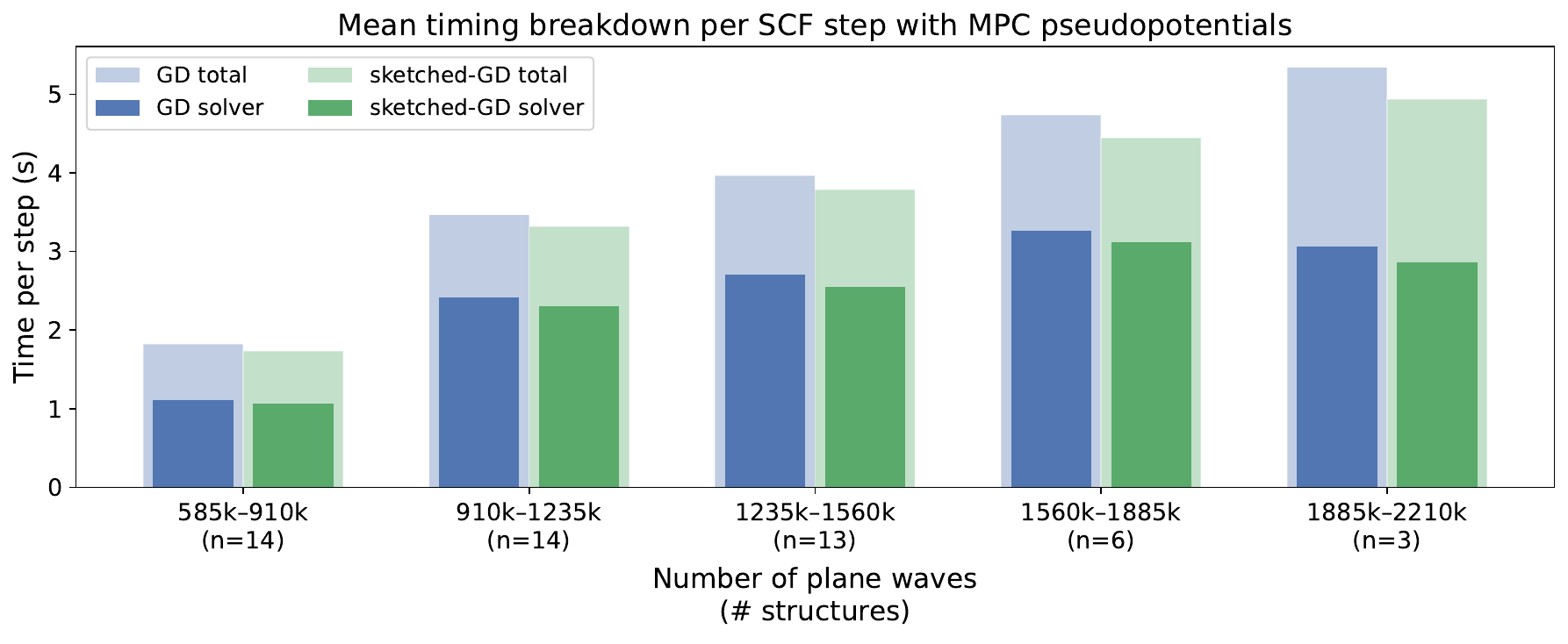}
    \caption{Mean GPU time per SCF step as a function of problem size in terms of number of plane waves. Lightly shaded bars represent total step time, while darker thinner bars represent the time spent in the eigensolver. Parenthesis on the x-axis labels show benchmark system distribution over the plane-wave bins.}
    \label{fig:dft_mpc_hist}
\end{figure*}

While there is a clear advantage to using the sketched generalized Davidson solver in this regime, it is interesting to breakdown the benchmark results by system size. Figure \ref{fig:dft_mpc_hist} displays the average time per SCF step across the 50 test structures, binned by problem size. As the number of plane waves increases, so does the cost of an SCF step. However, the randomized solver has a lower mean time than its deterministic counterpart in every reported size bin. We expect this trend to continue for larger systems, but limited GPU memory has restricted the scale of our benchmarks. For the three largest systems in this set, the randomized solver is 9\% faster. The cubic scaling of plane-wave DFT is another constraint, and we expect the method to work well with lower scaling basis choices such as wavelets or Gaussians. Finally, we note that using the randomized solver never degrades performance. In the worst case, there is simply no clear advantage. Over long SCF cycles or large high throughput studies, the compounding efficiency gains can amount to a significant amount of resources saved. 

\section{Conclusions}
\label{sec:Conclusions}

We have presented a sketched variant of the block generalized Davidson eigensolver, implemented it for CPUs and GPUs in mixed precision as part of the open-source \texttt{RandESC} package, and benchmarked it against the deterministic method it modifies and against LOBPCG, the incumbent solver of plane-wave DFT in DFTK. The benchmarks cover two main regimes: sparse test matrices, where the number of wanted eigenpairs $k$ is held fixed while the problem dimension $n$ grows, and self-consistent field DFT calculations with DFTK, where $k$ grows together with $n$. A third regime was investigated, in which DFT calculations were performed with a very high plane-wave cutoff, forcing $k \ll n$ while keeping $k$ and $n$ proportional.

The sketched variant replaces full orthogonalization by the randomized Gram-Schmidt process of Algorithm~\ref{alg:rgs}, which leaves the search space sketch-orthonormal rather than orthonormal. This in turn changes what the Rayleigh-Ritz step solves, from the standard eigenproblem \eqref{eq:rr} to the generalized Hermitian eigenproblem \eqref{eq:ghep} (Sec.~\ref{sec:extraction}), at a saving of $njk$ per iteration traded for a $j^3$ price in the dense solve (Sec.~\ref{sec:cost}). Extraction remains exact throughout: the resulting Ritz pairs are true Ritz pairs of $\vb A$, and the classical interlacing and convergence theory of Sec.~\ref{sec:JD} carries over unchanged. On the sparse matrices, with $k=500$ fixed and $n$ growing, the orthogonalization saving eventually dominates the added cost of the generalized eigenproblem: the sketched solver overtakes the deterministic one near $n \approx 2\cdot 10^4$ and is $25\%$ faster at $n=5\cdot 10^5$ (Sec.~\ref{sec:sparse_benchmark}).

In the first DFT benchmarks with DFTK (Sec.~\ref{sec:dft_benchmark}), both generalized Davidson variants beat the LOBPCG reference by 5 to 11\% in total SCF time. Sketching itself brings little further gain there: because the number of requested bands $k$ grows with the system size, as does $n$, the orthogonalization saving and the added cost of the generalized eigenproblem grow at the same rate and largely cancel, leaving the sketched and deterministic solvers within a few percent of each other.

In the second benchmark with DFTK (Sec.~\ref{sec:mpc_pseudos}), the MPC pseudopotentials are used. Near all-electron accuracy can be reached, but at the cost of a very high plane-wave cutoff. In this regime, $k$ and $n$ still grow proportionally, but $k \ll n$. We observe that both deterministic and sketched generalized Davidson implementation outperform the reference LOBPCG by 18\% and 19\% respectively. Additionally, we observe a clear advantage of the randomized solver, which grows with system size.

The regime in which sketching pays off is thus mainly set by how $k$ scales with $n$: a fixed, large $k$ against a growing $n$ favors sketching. However, in plane-wave DFT where $k \propto n$, the picture is more nuanced. With standard pseudopotentials and plane-wave cutoffs, the benefits of randomization are bounded by the $O(j^3)$ dense solve that grows at the same rate as the saving. However, when very high cutoffs are required and $k \ll n$, the sketched eigensolver becomes measurably favorable.

\begin{acknowledgments}
    We thank Raffaele Solc\`a, Simon Pintarelli, Anton Kozhevnikov, and Michael Herbst for fruitful discussions. 
    
We acknowledge support from the RandESC project (Novel algorithms based on randomization and mixed precision for electronic-structure calculations), via a grant by the Swiss National Supercomputing Centre (CSCS) as part of the Platform for Advanced Scientific Computing (PASC). This research was supported by the NCCR MARVEL, a National Centre of Competence in Research, funded by the Swiss National Science Foundation (grant number 205602), and by the Swiss National Science Foundation Grant No.~200021-227641 and No.~200021-236507. Computer time was provided by the Swiss National Supercomputing Centre (CSCS) under project No.~lp18, s1326, s1335, and mr33.
\end{acknowledgments}

\section*{AUTHOR DECLARATIONS}

\subsection*{Conflict of Interest}

The authors have no conflicts to disclose.









\section*{DATA AVAILABILITY}
The data that support the findings of this study are openly available in the RandESC repository~\cite{RandESC_repo}. The structures used for the GPU benchmarks are available at \url{https://github.com/RandESC-PASC/RandESC/tree/benchmarks/benchmark/mc3d_optimade}, the structures used for the CPU benchmarks at \url{https://github.com/RandESC-PASC/RandESC/tree/benchmarks/benchmark/mc3d\_optimade_cpu}, and the scripts and parameters used for the DFTK benchmarks are available at \url{https://github.com/RandESC-PASC/RandESC/tree/benchmarks/benchmark/DFTK}.

\bibliography{references}

\end{document}